\documentclass{amsart}
\usepackage{amssymb,amsrefs}

\usepackage[colorlinks]{hyperref}
\usepackage{booktabs}
\usepackage[all]{xy}
\usepackage{multirow}
\usepackage{booktabs,array,longtable}
\usepackage{tikz-cd}
\usepackage{mathtools}

\theoremstyle{plain}
\newtheorem{theorem}{Theorem}[section]

\newtheorem{lemma}[theorem]{Lemma}

\newtheorem{definition}[theorem]{Definition}

\theoremstyle{definition}
\newtheorem{remark}[theorem]{Remark}
\newtheorem{example}[theorem]{Example}

\DeclareMathOperator{\End}{End}
\DeclareMathOperator{\Frob}{Frob}
\DeclareMathOperator{\id}{id}
\DeclareMathOperator{\GL}{GL}

\DeclareMathOperator{\Gal}{Gal}

\newcommand{\FF}{{\mathbf{F}}}

\newcommand{\Que}{{\mathbf{Q}}}
\newcommand{\Zee}{{\mathbf{Z}}}

\newcommand{\gothp}{{\mathfrak{p}}}

\newcommand{\iso}{\cong}

\newcommand{\mapright}[1]{\mathop{\longrightarrow}\limits^{#1}}
\def\isar{\ \mapright{\sim}\ }

\renewcommand{\tilde}{\widetilde}

\newcommand\lowtilde{\lower0.7ex\hbox{\textasciitilde}}

\begin{document}

\title[Entanglement for elliptic curve group orders]
{Entanglement for elliptic curve group orders}

\author[Campagna]{Francesco Campagna}

\author[Jones]{Nathan Jones} 

\author[Pappalardi]{Francesco Pappalardi}

\author[Stevenhagen]{Peter Stevenhagen}

\address{Laboratoire de Math\'ematiques Blaise Pascal,
         Universit\'e Clermont-Auvergne,
         Cler\-mont-Ferrand,
         France}
         \email{francesco.campagna@uca.fr}
         
\address{Department of Mathematics, Statistics and Computer Science,
         University of Illinois at Chicago, 
         851 S Morgan St, 322 SEO,
         Chicago, IL 60607,
         USA}
         \email{ncjones@uic.edu}
        
\address{Dipartimento di Architettura,
         Universit\`a Roma Tre,
         Largo Giovanni Battista Marzi 10,
         I--00156 Rome,
         Italy}
         \email{francesco.pappalardi@uniroma3.it}
         
\address{Mathematisch Instituut,
         Universiteit Leiden, 
         Postbus 9512, 
         2300 RA Leiden, The Netherlands}
         
\email{psh@math.leidenuniv.nl}

\date{\today}


\begin{abstract}
For an elliptic curve $E$ defined over a number field $K$ and
$\gothp$ a prime of good reduction of $E/K$, we let 
$n_\gothp=\#E(k_\gothp)$ be the order of the group of points over 
the residue class field $k_\gothp$.
Divisibility of $n_\gothp$ by a prime $\ell$ depends on 
the splitting behavior of $\gothp$ in the $\ell$-division field 
$K(E[\ell](\overline K)$) of $E$ and, 
generically, divisibility of $n_\gothp$ by different primes $\ell_1$ 
and $\ell_2$ are `independent events' on the set of primes $\gothp$ 
of good reduction.
However, there do exist elliptic curves with `entangled division fields' 
for which infinitely many $n_\gothp$ are coprime to $\ell_1$ and 
infinitely many $n_\gothp$ are coprime to $\ell_2$,
but {\it no\/} $n_\gothp$ is coprime to to $\ell_1 \ell_2$.
For $\ell_1=2$, we classify in group theoretical terms the entanglement 
giving rise to this behavior, and count up to isomorphism the number 
of `special $2\ell$-division groups' exhibiting the phenomenon.
A final section shows how such groups may be realized as the
Galois group of $2\ell$ division fields for elliptic
curves $E/\Que$.

\end{abstract}

\maketitle

\section{Introduction}
\label{S:intro}

\noindent
Let $K$ be a number field and $\gothp$ a prime of $K$ of norm $N\gothp$.
Then the unit group $k_\gothp^*$ of 
the residue class field at $\gothp$ is a cyclic group of order $m_\gothp=N\gothp-1$.
A prime number~$\ell$ divides $m_\gothp$ if and only if $\gothp$ is a prime of $K$ not dividing $\ell$ that splits completely in the $\ell$-th cyclotomic extension $K\subset K_\ell=K(\zeta_\ell)$.
This formulation implies that the set of primes $\gothp$ with $\ell|m_\gothp$
has a positive density $[K(\zeta_\ell):K]^{-1}$ inside the set of 
all primes of $K$. 
In particular, almost all groups $k_\gothp^*$ contain an element of order $\ell$ if and only if $K^*$ contains an element of order $\ell$.

Now let $E$ be an elliptic curve defined over $K$, and  
$\gothp$ a prime of good reduction of $E$.
Then the divisibility by $\ell$ of the order $n_\gothp=\#E(k_\gothp)$ 
of the group of points over $k_\gothp$ depends
on the splitting behavior of $\gothp$ in the extension $K\subset K_\ell$
obtained by taking $K_\ell=K(E[\ell](\overline K))$ to be the
$\ell$-division field of $E$.
This is the field obtained by adjoining to $K$ 
the coordinates of all the $\ell$-torsion points of $E$ defined over an 
algebraic closure $\overline K$ of $K$.

In the multiplicative case, the condition for the divisibility
$\ell|m_\gothp$ at a prime $\gothp\nmid \ell$ simply
amounts to requiring that the Frobenius of $\gothp$ in the 
cyclotomic Galois group $\Gal(K_\ell/K)$ be trivial:
\begin{equation}\label{mult-div}
\ell\mid m_\gothp \Leftrightarrow \Frob_\gothp={\rm id}.
\end{equation}
In the elliptic case, the condition for the divisibility
$\ell|n_\gothp$ at a prime $\gothp\nmid \ell$ of good reduction
amounts to requiring that the Frobenius of $\gothp$ in the 
$\ell$-division Galois group $G_\ell=\Gal(K_\ell/K)$, 
which is well-defined up to conjugacy, leave invariant some non-zero 
$\ell$-torsion point in $V_\ell=E[\ell](\overline K)$.
Define, for $\sigma\in G_\ell$, the \emph{fixpoint rank} 
\begin{equation}
\label{eqn:fprank-l}
r(\sigma)= \dim_{\FF_\ell} V_\ell^{\langle\sigma\rangle}
= \dim_{\FF_\ell} E[\ell](\overline K)^{\langle\sigma\rangle}
\in \{0, 1, 2\}
\end{equation}
of $\sigma$ as the dimension of the subspace of $\sigma$-invariant $\ell$-torsion points.
Then for primes $\gothp\nmid \ell$ of good reduction,
the analogue of \eqref{mult-div} is
\begin{equation}
\label{ell-div}
\ell\mid n_\gothp \Leftrightarrow r(\Frob_\gothp)\ge1.
\end{equation}
It is a `weaker' condition than \eqref{mult-div}, as it states 
that $\Frob_\gothp$ 
fixes \emph{some} point of order $\ell$ rather than all of them.

All conjugacy classes in $G_\ell$ occur as the Frobenius class of infinitely many primes,
so an elliptic curve $E/K$ has infinitely many primes $\gothp$ of good reduction
with group order $n_\gothp=\#E(k_\gothp)$ coprime to $\ell$ 
if and only if $G_\ell$ contains elements of fixed point rank 0.
If $G_\ell$ is without such elements 
and almost all $E(k_\gothp)$ contain an element of order $\ell$, 
this does \emph{not} in general imply that $E(K)$ contains an element of order $\ell$.
There is also the possibility that $E$ is $\ell$-isogenous
over $K$ to an elliptic curve $E'/K$ for which
$E(K)$ contains an element of order $\ell$, cf.\ Remark \ref{serre-ex-remark}.
\bigskip

Now suppose $\ell_1$ and $\ell_2$ are different prime numbers.
If there exist primes $\gothp_1, \gothp_2$ of $K$ 
not dividing $\ell_1\ell_2$ such that the group order 
$m_{\gothp_i}=\#k_{\gothp_i}^*$ is coprime to $\ell_i$ for $i=1, 2$,
we can find infinitely many primes $\gothp$ in $K$ for which 
$m_\gothp$ is coprime to $\ell_1\ell_2$.
Using \eqref{mult-div}, this follows from the trivial 
group theoretical fact that the Galois group 
$G_{\ell_1\ell_2}=\Gal(K(\zeta_{\ell_1\ell_2})/K)$, when viewed as a subgroup of the product
\[
G_{\ell_1\ell_2} \subset G_{\ell_1}\times G_{\ell_2}=\Gal(K(\zeta_{\ell_1})/K)\times \Gal(K(\zeta_{\ell_2})/K),
\]
contains an element with two non-identity components if and only if it contains 
an element with a non-identity first component and an element with a non-identity second component.

For elliptic group orders, the non-triviality of Frobenius elements in
$G_\ell\subset \GL(V_\ell)$ is replaced by having fixpoint rank zero, and 
the situation is very different.
In this case the Galois group $G_{\ell_1\ell_2}$
of the $\ell_1\ell_2$-division field of $E$ over $K$ embeds as a subgroup 
\[
G_{\ell_1\ell_2} \subset G_{\ell_1}\times G_{\ell_2}
\subset \GL(V_{\ell_1}) \times \GL(V_{\ell_2})
\iso \GL_2(\FF_{\ell_1}) \times \GL_2(\FF_{\ell_2})
\]
in the product of
the Galois groups of the $\ell_1$- and $\ell_2$-division fields of $E$ 
over~$K$, with 
$G_{\ell_i}={\rm im}[G_{\ell_1\ell_2}\to \GL(V_{\ell_i})]$ 
acting linearly on the 2-dimensional 
$\FF_{\ell_i}$-vector space $V_{\ell_i}$ of $\ell_i$-torsion elements.
Note that any subgroup $G$ of the group on the right comes
with two \emph{fixpoint rank functions}
$r_{\ell_i}: G\to \{0, 1, 2\}$ that are obtained by composing the projection on the $\ell_i$-component with the map 
\eqref{eqn:fprank-l} for $\ell=\ell_i$. 
We can define purely group theoretically when a 
group of this nature deserves to be called \emph{special}.

\begin{definition} \label{def:special_group}
Let $V_{\ell_1}$ and $V_{\ell_2}$ be two-dimensional vector spaces over
$\FF_{\ell_1}$ and $\FF_{\ell_2}$.\break
    A subgroup 
    \[
    G\subset \GL(V_{\ell_1}) \times \GL(V_{\ell_2})
    \]
    is called {\bf special of level $\ell_1\ell_2$\/} if both fixpoint rank functions
    $r_{\ell_i}: G\to \{0, 1, 2\}$ assume the value $0$ on $G$, but
    not simultaneously:
\[
(r_{\ell_1}(g), r_{\ell_2}(g))\ne(0, 0) 
\quad\text{ for all~$g \in G$}.
\]
\end{definition}
If $G$ arises as the Galois group of the $\ell_1\ell_2$-division field 
of an elliptic curve $E$ over a number field $K$, 
the implication of $G$ being special is that the group orders 
$n_\gothp=\#E(k_\gothp)$ are infinitely often 
coprime to $\ell_1$ and to $\ell_2$,  but not to $\ell_1\ell_2$.
This is the motivation for us to abstractly classify such groups, and then to realize them as Galois groups of $\ell_1\ell_2$-division fields.

\medskip\noindent
Writing $G_{\ell_i}$ for the projection of $G$ on the $\ell_i$-component,
the \emph{Goursat lemma} from group theory tells us that we may describe 
groups $G$ as in Definition \ref{def:special_group}
as the graph of an isomorphism of quotients
\begin{equation}
\label{eqn:Goursat}
     \psi: G_{\ell_1}/N_{\ell_1} \xrightarrow{\sim} G_{\ell_2}/N_{\ell_2}.
\end{equation}
In more explicit terms: $G$ is the subgroup
\begin{equation}\label{eqn:Goursatform}
G=\{(\sigma_1, \sigma_2): \psi (\sigma_1\bmod N_{\ell_1})
   = \sigma_2\bmod N_{\ell_2}\}
   \subset G_{\ell_1}\times G_{\ell_2}.
\end{equation}   
We have an index-equality
\[
[G_{\ell_1}\times G_{\ell_2} : G]=
[G_{\ell_i}:N_{\ell_i}].
\]
For special $G$, the inclusion 
$G\subset G_{\ell_1}\times G_{\ell_2}$ is necessarily strict.

In the case of level $\ell_1\ell_2=2\ell$, with $\ell$ an odd prime,
we can classify the special groups.
It turns out that for $G\subset G_2 \times G_\ell$
special, the component group $G_2$ is the full group
$\GL(V_2)$, whereas $G_\ell$ is contained in a Borel subgroup
$B\subset \GL(V_\ell)$.
The Goursat isomorphism \eqref{eqn:Goursat} describing 
$G$ is an isomorphism between groups of order~2, and can be described 
explicitly in terms of the unique quadratic character 
\[
\varepsilon: \GL(V_2) \cong S_3 \to \mathbb{F}_\ell^\ast
\]
and the fundamental characters
$\chi_1,\chi_2: B \to \mathbb{F}_\ell^\ast$ of
a Borel group containing $G_\ell$.
(See the next section for the definition of 
these fundamental characters $\chi_i$.)
\begin{theorem} 
\label{thm:2l_entanglement}
Suppose $G\subset \GL(V_2) \times \GL(V_\ell)$ 
is a special subgroup of level $2\ell$.
Then $G_2$ equals $\GL(V_2)$, and $G_\ell$ is contained in
a Borel subgroup $B\subset \GL(V_\ell)$.

The restrictions 
$\chi_1,\chi_2: G_\ell \to \mathbb{F}_\ell^\ast$
of the fundamental characters of $B$ to $G_\ell$ are
non-trivial, and
there exists $i\in\{1,2\}$ such that $\chi_i$
is quadratic on $G_\ell$ and
\begin{equation}\label{eqn:def-of-G}
    G=\{(\sigma_2,\sigma_\ell) \in \GL(V_2) \times G_\ell: \varepsilon(\sigma_2)=\chi_i(\sigma_\ell)\}.
\end{equation}
Conversely, for a subgroup $G_\ell\subset B$ for which the fundamental characters are both non-trivial, with $\chi_i$ quadratic,
equation \eqref{eqn:def-of-G} defines a special subgroup of level~$2\ell$.
\end{theorem}
\noindent
Theorem \ref{thm:2l_entanglement} is so explicit that we can actually count the number of such groups, taken up to conjugation.
\begin{theorem} 
\label{thm:countingspecials}

Up to conjugation in $\GL(V_2) \times \GL(V_\ell)$, 
the number of special subgroups of level $2\ell$ equals 
\[
\textstyle
3\,\sigma_0(\ell-1)+3\, \sigma_0\left(\frac{\ell-1}{2}\right)-4.
\]
Here $\sigma_0(k)$ denotes the number of positive divisors of $k$.

\end{theorem}
\noindent
Section \ref{S:Borel} collects basic facts on 
subgroups of $\GL(V)$ for a 2-dimensional vector space~$V$, including a Lemma \ref{ai-lemma} on fixpoint ranks 
that was suggested to us by AI.
It replaces our initial longer argument that was based on modular representation theory.
The proof of Theorem \ref{thm:2l_entanglement} classifying 
special subgroups of level $2\ell$ is given in 
Section \ref{S:mainproof}, and Section \ref{S:countingformula} provides 
a counting formula for them.
A final section \ref{S:ellipticcurves}
discusses their realisations by elliptic curves defined over~$\bf Q$.

\bigskip\noindent
{\bf Acknowledgements.}
Francesco Pappalardi and Peter Stevenhagen received support from the grants \textit{AAP Accueil de conférenciers invités internationaux} and \textit{AAP Impulsion} of Universit\'e Clermont Auvergne (UCA) while working on this paper.
They thank the UCA for its financial support and the Laboratoire de Math\'ematiques Blaise Pascal for its inspiring atmosphere.

\section{Linear group actions in dimension 2}
\label{S:Borel}

Let $G$ be a group acting linearly on a finite dimensional $k$-vector space $V$.
Then the \emph{fixpoint rank} $r:G \to\Zee$ associated to this action is 
the function defined by
\begin{equation}
\label{eqn:fprank}
r(g)= \dim_k V^{\langle g\rangle}.
\end{equation}
It is bounded by $\dim_k(V)$, and in this section we will 
restrict to the special case where the action 
$\rho: G\to \GL(V)$ is a 2-dimensional representation: $\dim_k(V)=2$.

In view of our criterion \eqref{ell-div}, we are interested in elements in $G$ 
having $r(g)\ge 1$.
The case where all elements $g\in G$ have $r(g)\ge 1$ is well known, and occurs already 
as an exercise in Serre's book \cite{Serre2}*{p.\ I-2}.
We formulate this result as Lemma \ref{serre-ex} below, and note that it has a useful 
generalisation to $n$-dimensional representations of groups $G$ for which
all elements $g\in G$ have $r(g)\ge n-1$, see \cite{JPapStev}*{Lemma 4.1}.
\begin{lemma}\label{serre-ex}
Suppose all elements $g\in G$ have positive fixpoint rank.
Then there exists a 1-dimensional subspace $L\subset V$ with the property that at least 
one of the following holds:
\begin{enumerate}
\item
all $g\in G$ act as the identity on $L$;
\item
$G$ stabilizes $L$, and all $g\in G$ act as the identity on $V/L$; \qed
\end{enumerate}
\end{lemma}
\noindent
\begin{remark}\label{serre-ex-remark}
For an elliptic curve $E$ over a number field $K$, 
the $\ell$-torsion subgroup $V=E[\ell](\overline K)$ is a 
2-dimensional vector space over the field of $\ell$ elements, and
the Galois group $G$ of the $\ell$-division field of $E$ over $K$ acts
linearly on $V$.
If almost all local point groups are of order
$n_\gothp=\# E(k_\gothp)$ divisible by $\ell$, then 
\eqref{ell-div} implies that the elements $g\in G$, which all occur as 
Frobenius elements of infinitely many~$\gothp$, 
have positive fixpoint rank.
Now Lemma \ref{serre-ex} implies that either $E$
itself has a $K$-rational point of order $\ell$, or $E$ is $\ell$-isogenous over $K$ to such a curve.
\end{remark}
Lemma \ref{serre-ex} describes $\rho[G]$ as a subgroup of 
the Borel group 
$B_L\subset \GL(V)$, which is defined as the stabilizer 
$B_L=\{\sigma\in\GL(V): \sigma[L]=L\}$ of the subspace $L$.
We have an exact sequence of $k[B_L]$-modules 
\[
0\to L \to V \to V/L \to 0
\]
in which the 1-dimensional representations $L$ and $V/L$
of $B_L$ are given by the \emph{fundamental characters} 
$\chi_1=\chi_L$ and $\chi_2=\chi_{V/L}$ mapping  $B_L$ to $k^*$.
The product $\chi_1\chi_2: B_L\to k^*$
is the restriction of the determinant map $\det: \GL(V)\to k^*$ to $B_L$.

Non-trivial elements $\sigma\in B_L$ with $\chi_1(\sigma)=\chi_2(\sigma)=1$
are called \emph{transvections}.
Together with the identity, they form a normal subgroup 
$T_L\subset \GL(V)$ isomorphic to the additive group $k$.
We have a split exact sequence of groups
\begin{equation}\label{Borel}
0\to T_L \longrightarrow B_L \mapright{(\chi_1, \chi_2)} k^* \times k^* \to 1,
\end{equation}
showing that $B_L$ is a semi-direct product of $T_L$ and a 
split Cartan subgroup in $\GL(V)$, representable in matrix form as
\[
B_L = 
\begin{pmatrix}
1 & k  \\
0 & 1  
\end{pmatrix}
\rtimes 
\begin{pmatrix}
k^* & 0 \\
0 & k^*  
\end{pmatrix}.
\]

In terms of the representation $\rho: G\to \GL(V)$,
Lemma \ref{serre-ex} states that $G$ acts on~$V$ 
by elements of positive fixpoint rank if and only if we have
\begin{equation}
  \rho[G]\subset K_i=\ker\chi_i\subset B_L  
\end{equation}
for a Borel subgroup $B_L$ and a fundamental character
$\chi_i$ on $B_L$.
In the case $\rho[G]\subset T_L$ the inclusion 
$\rho[G]\subset K_i$ holds for both 
fundamental characters on $B_L$.

A Borel subgroup $B_L$ is it own normalizer in $\GL(V)$, 
its conjugates in $\GL(V)$
being of the form $\tau B_L\tau^{-1}=B_{\tau L}$.
In the same way, one sees that the
kernels $K_i=\ker\chi_i\subset B_L$ of the fundamental characters 
have $B_L$ as their normalizer in $\GL(V)$.
For subgroups of $K_i$ we have the following result.
\begin{lemma}\label{borel-normalizer}
For $i\in\{1,2\}$, the normalizer in $\GL(V)$ of a non-trivial subgroup 
$H$ of $K_i=\ker\chi_i\subset B_L$ is contained in $B_L$.
\end{lemma}
\begin{proof}

The proof for $K_1=\ker\chi_1$ may be given by staring at the matrix identity
\[
\begin{pmatrix}
a & b  \\
c & d  
\end{pmatrix}
\begin{pmatrix}
1 & x  \\
0 & y  
\end{pmatrix}
=
\begin{pmatrix}
1 & \tilde x  \\
0 & \tilde y  
\end{pmatrix}
\begin{pmatrix}
a & b  \\
c & d  
\end{pmatrix},  
\]
which yields identities $a=a+c\tilde x$ and $c=c\tilde y$,
whence $c\tilde x = c (1-\tilde y)=0$ and $c=0$.
Transposing all matrices, the statement for $K_2=\ker\chi_2$ is obtained. 
\end{proof}
\noindent
The fixpoint rank is not a function on a group with nice algebraic properties,
but there are ways to produce new elements of positive fixpoint rank from old ones.
The following elementary lemma (and its elementary proof!) 
were obtained with the help of ChatGPT.
It replaces our original arguments based on modular representation theory.
This is the only instance of AI used in this paper.
\begin{lemma}\label{ai-lemma}
   For elements $g, x \in \GL(V)$, suppose that $g$, $gx$ and $g^2x$ 
   have positive fixpoint rank.
   Then $x$ has positive fixpoint rank. 
\end{lemma}
\begin{proof}
    The characteristic polynomial of an element $A\in \GL(V)$ can be written as 
    \[
    f_A=\det(A-X\cdot\id)=(X-\lambda_1)(X-\lambda_2)=X^2-t_A X+d_A \in k[X],
    \]
    with $t_A$ and $d_A$ the trace and the determinant of $A$, and 
    $\lambda_1$, $\lambda_2\in\overline k$ its eigenvalues.
    We have $r(A)\ge1$ if and only if $A$ has at least one eigenvalue 1, 
    which amounts to having $f_A(1)=0$.
    Hence 
    \begin{equation}\label{eqn:d=t-1}
    r(A)\ge1 \Longleftrightarrow 
    t_A=1+d_A.
    \end{equation}
By the Cayley-Hamilton identity, multiplying the zero element 
$f_g(g)=0\in\End_k(V)$ by $x$ yields
\[
g^2x-t_g gx+d_gx=0\in \End_k(V).
\]
Taking traces, we have
\[
t_{g^2x}-t_gt_{gx}+d_gt_x=0\in k.
\]
As $g$, $gx$ and $g^2x$ satisfy \eqref{eqn:d=t-1} by assumption, 
we can rewrite their traces $t$ as $1+d$ and use the multiplicativity of the 
determinant $d$ to miraculously obtain
\[
1+d_{g^2x}-(1+d_g)(1+d_{gx})+d_gt_x=
-d_g(1+d_x-t_x)=0\in k
\]
Dividing by $-d_g\in k^*$, we see that both sides of \eqref{eqn:d=t-1} are satisfied for $A=x$.
\end{proof}

\noindent
We now specialize to the case where $V=V_\ell$ is of dimension 2
over the finite field $k=\FF_\ell$ of prime order $\ell$.
In this case $\GL(V_\ell)$ is a finite group of order 
\[
\#\GL(V_\ell)=\ell(\ell+1)(\ell-1)^2.
\]
For each of the $\ell+1$ lines $L\subset V_\ell$ through 0, the
corresponding Borel subgroup $B_L\subset \GL(V_\ell)$ of index $\ell+1$
has order $\ell(\ell-1)^2$.
Any transvection in $B_L$ generates the normal $\ell$-Sylow subgroup $T_L\subset B_L$ 
of order $\ell$.
The $(\ell+1)(\ell-1)$ transvections in $\GL(V_\ell)$ form a single conjugacy 
class in $\GL(V_\ell)$ consisting of the elements in $G$ of order $\ell$.

\section{Special groups of level 2\texorpdfstring{$\ell$}{ell}}
\label{S:mainproof}
\noindent
For the proof of Theorem \ref{thm:2l_entanglement}, we let
the subgroup 
    \[
    G\subset \GL(V_2) \times \GL(V_\ell)
    \]
be special in the sense of Definition \ref{def:special_group}.
Denote the projections of $G$ on its two components
by $G_2$ and $G_\ell$, and write the Goursat isomorphism defining $G$ as
in \eqref{eqn:Goursatform} as
\begin{equation}\label{Goursat-psi}
     \psi:   G_\ell/N_\ell \xrightarrow{\sim} G_2/N_2.
\end{equation}
As $G$ is special, this is an isomorphism between non-trivial groups.

\medskip\noindent
{\bf Step 1}: we have $G_2=\GL(V_2)\iso S_3$.

\smallskip\noindent
{\bf Proof}.
The elements in $\GL(V_2)$ of fixpoint rank 0 are those that 
permute the 3 non-zero elements of $V_2$ in a 3-cycle, so $G_2$ contains
the subgroup $A_3\subset \GL(V_2)$ of order~3.
If $G_2$ is not the full group $\GL(V_2)\iso S_3$, the 
Goursat isomorphism \eqref{Goursat-psi} is 
\[
\psi: G_\ell/N_\ell\isar G_2=A_3. 
\]
Pick $g\in G_\ell\setminus N_\ell$ and $x\in N_\ell$.
As $\psi$ maps the residue classes of the elements $g$, $gx$ and $g^2x$
modulo $N_\ell$ to elements of fixpoint rank 0, 
these elements have positive fixpoint rank by the assumption that $G$ be special.
Now Lemma \ref{ai-lemma} implies that $x$ also has positive fixpoint rank, and we conclude that all elements in $G_\ell$ have positive fixpoint rank. 
Contradiction.\qed

\medskip\noindent
{\bf Step 2}: $\psi: G_\ell/N_\ell\isar G_2/N_2$ is an isomorphism between groups of order 2.

\smallskip\noindent
{\bf Proof}.
As we have $G_2\iso S_3$, its only non-trivial quotient of 
order not 2 is $G_2$ itself. 
In order to derive a contradiction, suppose we have a Goursat isomorphism 
\[
\psi: G_\ell/N_\ell\isar G_2= S_3,
\]
and write $\psi^{-1}[A_3]=A/N_\ell$.
By the argument in Step 1, all elements $x\in A$ have positive fixpoint rank.
By Lemma \ref{serre-ex}, we find 
$A\subset K_i=\ker\chi_i\subset B_L$, 
so Lemma~\ref{borel-normalizer} implies $G_\ell\subset B_L$.

Suppose first that we have $\ell>3$.
Then all transvections in $G_\ell$ are trivial in the quotient $G_\ell/N_\ell$ of order 6, so they are contained in $N_\ell$.
As $B_L/T_L$ is abelian, this implies that $G_\ell/N_\ell$ 
is abelian, and that 
$\psi[G_\ell/N_\ell]$ is not $S_3$: 
contradiction.

In the remaining case $\ell=3$, the group $G_3$ 
of order divisible by 3 contains $T_L$ as its 3-Sylow subgroup. 
By the previous argument, $T_L$ is not contained in $N_3$, 
so we have $A=T_L N_3$.
The conjugation action of any element $g\in G_3$ outside $A$ on 
$A/N_3$ is by inversion, so $g$ is not in 
${\bf F}_3^* \cdot T_L$ and its eigenvalues in ${\bf F}_3^*$
are 1 and $-1$.
This implies $r(g)\ge1$, so $G_3$ only has elements 
of positive fixpoint rank: contradiction.

We conclude that the Goursat isomorphism \eqref{Goursat-psi} for $G$ 
is between groups of order~2.
We can view it as a quadratic character $\psi: G_\ell\to\FF_\ell^*$
that describes $G$ as 
\begin{equation}
\label{eqn:G-description}
    G=\{(\sigma_2,\sigma_\ell) \in \GL(V_2) \times G_\ell: 
    \psi(\sigma_\ell)=\varepsilon(\sigma_2)\},
\end{equation}
with $\varepsilon: \GL(V_2) \cong S_3 \to \mathbb{F}_\ell^\ast$ 
the unique quadratic character on $\GL(V_2)$. \qed

\medskip\noindent
{\bf Step 3}: we have $G_\ell\subset B_L$ for a Borel subgroup in $\GL(V_\ell)$, the restrictions of the fundamental characters $\chi_1$ or $\chi_2$ to $G_\ell$ are non-trivial, and one of them
equals
the quadratic character $\psi$ in \eqref{eqn:G-description}.

\smallskip\noindent
{\bf Proof}.
As elements $\sigma_\ell\in N_\ell=\ker\psi$ can be paired with 
the elements  $\sigma_2\in \ker\varepsilon$ of fixpoint rank 0 
to obtain $(\sigma_2,\sigma_\ell)\in G$,
we find that all elements of $N_\ell$ have fixpoint rank at least 1.
By Lemma \ref{serre-ex}, 
this implies that $N_\ell$ is contained in a
subgroup $K_i=\ker\chi_i$ of some Borel group $B_L$.

If $N_\ell=1$ is the trivial group, then $G_\ell$ 
is generated by an element $\sigma_\ell$ of order~2 
of fixpoint rank 0. 
As $\sigma_\ell$ is without eigenvalue 1 and $\ell$ is odd, 
we have $\sigma_\ell=-\id$.
This is a case of Theorem \ref{thm:2l_entanglement} with
$\psi=\chi_1=\chi_2: G_\ell\to\FF_\ell^*$ non-trivial quadratic. 
Here $G_\ell$ is contained in {\it all} Borel subgroups of $\GL(V_\ell)$.

If $N_\ell$ is non-trivial and contained in a
subgroup $K_i=\ker\chi_i$ of some Borel group $B_L$,
then Lemma \ref{borel-normalizer}
implies that we have $G_\ell\subset B_L$.

In the case where $N_\ell\subset K_1=\ker\chi_1\subset B_L$ fixes
$L$ pointwise,
there exists $g\in G_\ell\subset B_L$ with $g^2\in N_\ell$ and $r(g)=0$. 
This implies that $g$ acts on $L$ with eigenvalue $-1$, 
as do all other elements in the coset 
$gN_\ell=G_\ell\setminus N_\ell$.
This shows that the restriction of $\chi_1$ to $G_\ell$ is a quadratic character that we may identify with $\psi$.
As elements $g\in G_\ell$ with $r(g)=0$ satisfy $\chi_2(g)\ne 1$,
the other character $\chi_2$ is non-trivial on $G_\ell$.

In the same way, the case $N_\ell\subset K_2=\ker\chi_2\subset B_L$ 
in which $N_\ell$ acts trivially on $V/L$ gives rise to 
$\chi_2$ being the quadratic character on $G_\ell\subset B_L$
that we may identify with $\psi$.
As elements in $G_\ell$ with $r(g)=0$ do not fix $L$ pointwise,
$\chi_1$ is non-trivial on $G_\ell$. \qed

\medskip\noindent
Step 3 shows that special subgroups of level $2\ell$ satisfy
the description given in Theorem \ref{thm:2l_entanglement}.
For the converse, define $G$ as in \eqref{eqn:G-description}
with $G_\ell\subset B_L$ and $\psi: G_\ell\to {\bf F}_\ell^*$
a non-trivial quadratic character obtained as the restriction
of a fundamental character $\chi_i$ on $B_L$,
and assume the other fundamental character $\chi_{3-i}$
is non-trivial on $G_\ell$.
Then all elements in $g\in N_\ell=\ker\psi\subset G_\ell$ have
an eigenvalue $\chi_i(g)=1$ and $r_\ell(g)\ge 1$.
As $\chi_{3-i}$ is non-trivial on $G_\ell$, it cannot be trivial
on $G_\ell\setminus N_\ell$, so there exist elements $g\in G_\ell$
with $\chi_i(g)=-1$ and $\chi_{3-i}(g)\ne 1$, whence $r_\ell(g)=0$.
In \eqref{eqn:G-description}, these elements are paired in $G$ 
with elements $g_2\in G_2$ that have $\varepsilon(g_2)=-1$
and $r_2(g)\ne 0$ as they fix a non-zero vector in $V_2$.
This shows that $G$ is special, and concludes the proof
of Theorem \ref{thm:2l_entanglement}. \qed
\section{Counting special subgroups of level \texorpdfstring{$2\ell$}{2l}}
\label{S:countingformula}
\noindent
The explicit nature of the classification provided by 
Theorem \ref{thm:2l_entanglement} allows us to count the number
of special subgroups $G$ of level $2\ell$, up to conjugation
in $\GL(V_2)\times \GL(V_\ell)$.
The description in the Theorem provided by 
equation \eqref{eqn:def-of-G} shows that every special subgroup
of level $2\ell$ is obtained in the following way:
\begin{enumerate}
    \item [(a)]
    pick a subgroup $G_\ell$ of a Borel subgroup $B\subset \GL(V)$ 
such that both fundamental characters $\chi_1$ and $\chi_2$ are 
non-trivial, and at least one of them is quadratic;
    \item [(b)]
    define $G$ as in \eqref{eqn:G-description}, with $\psi$ 
    equal to a quadratic fundamental character.
\end{enumerate}
As the pair $(G_\ell, \ker\psi)$ completely determines the
special subgroup $G$, the counting can take place inside 
the Borel subgroup $B$.
For $B$, we have a surjection $\pi=(\chi_1,\chi_2): B\to \FF_\ell^*\times \FF_\ell^*$
giving rise to the split exact sequence \eqref{Borel} describing
$B$ as the semi-direct product of a subgroup
$\FF_\ell^*\times \FF_\ell^*$ of `diagonal matrices' 
with a normal transvection subgroup $T\iso\FF_\ell$ in $B$.
\[
\begin{tikzcd}[column sep=large, row sep=large]
1 \arrow[r]
  & T \arrow[r]
  & B \arrow[r, "\pi", "{(\chi_1,\chi_2)}"']
  & \mathbb{F}_\ell^{*} \times \mathbb{F}_\ell^{*} \arrow[r]
  & 1 \\
1 \arrow[r]
  & T_\ell \arrow[r] \arrow[u, phantom, sloped, "\subseteq"]
  & G_\ell \arrow[r, "\pi"] \arrow[u, phantom, sloped, "\subseteq"]
  & \pi[G_\ell] \arrow[r] \arrow[u, phantom, sloped, "\subseteq"]    
  & 1 \\
  &
  & N_\ell \arrow[r]
      \arrow[u, phantom, sloped, "\subseteq"]
      \arrow[u, phantom, "\scriptstyle 2"{xshift=2ex, yshift=-0.1ex}]
  & \pi[N_\ell]
      \arrow[u, phantom, sloped, "\subseteq"]
      \arrow[u, phantom, "\scriptstyle 2"{xshift=2ex, yshift=-0.1ex}]
  &
\end{tikzcd}
\]
In a similar way, we obtain a description of $G_\ell$ as the semi-direct product
of a subgroup $\pi[G_\ell]\subset \FF_\ell^*\times \FF_\ell^*$ of diagonal matrices
with the transvection subgroup $T_\ell=G_\ell\cap T$.

Given $G_\ell$, conjugation in $B$ fixes the normal subgroup $T$, 
and it does not change the image $\pi[G_\ell]$ in the abelian group $\FF_\ell^*\times \FF_\ell^*$.
In fact, the group $T_\ell$, which is of order either 1 or $\ell$, and the
group $\pi[G_\ell]$ together determine $G_\ell$ up to conjugation
inside $B$. 
Indeed, suppose the subgroups $G_\ell$ and $G_\ell'$ of $B$ 
satisfy $\pi[G_\ell]=\pi[G_\ell']$
and $G_\ell\cap T=G_\ell' \cap T=T_\ell$. 
If $T_\ell$ is trivial, 
then $G_\ell$ and $G_\ell'$ are both complements of $T$ in $G=G_\ell T_L=G_\ell' T_L$. 
As $G_\ell$ and $G_\ell'$ have order coprime to $\#T=\ell$, 
the Schur-Zassenhaus theorem \cite{Rotman}*{Theorem 7.40} 
implies that they are conjugate by an element of $G \subset B$. 
If $T_\ell$ is of order $\ell$, then $\pi[G_\ell]=\pi[G_\ell']$ actually implies $G_\ell=G_\ell'$.

As every special group $G$ of level $2\ell$ is constructed from
a subgroup $G_\ell$ of a Borel subgroup in $\GL(V_\ell)$
and a subgroup $N_\ell=\ker\psi\subset G_\ell$ of index 2,
we will approach the counting of special subgroups  
by going over the subgroups 
\begin{equation}
    \label{eqn:Goursatinclusion}
    \pi[G_\ell] \subset \chi_1[G_\ell]\times \chi_2[G_\ell]
\end{equation}
inside $\FF_\ell^*\times \FF_\ell^*$
meeting the conditions of Theorem \ref{thm:2l_entanglement}
and determining how many pairwise non-conjugate pairs 
$(G_\ell, N_\ell)$ we obtain in each case. 
We will have to deal with the `Goursat-index' 
\[
I(G)=[ \chi_1[G_\ell]\times \chi_2[G_\ell] : \pi[G_\ell]]
\]
of the inclusion \eqref{eqn:Goursatinclusion} and the difference between $G_\ell$ and 
$\pi[G_\ell]=G_\ell/T_\ell$.

\medskip\noindent
{\bf Proof of Theorem \ref{thm:countingspecials}.}
Let $G_\ell\subset B$ be a subgroup on which both fundamental characters 
$\chi_1$ and $\chi_2$ are non-trivial, and one of them
is quadratic with kernel $N_\ell$.
This gives rise to a special subgroup of level $2\ell$ as in Theorem \ref{thm:2l_entanglement},
with $N_\ell=\ker\chi_i$ the `entanglement kernel' at $\ell$ from \eqref{eqn:Goursat}.
For $m$ dividing $\ell-1$,
denote by $C_m\subset {\bf F}_\ell^*$ the subgroup of order 
$m$ in~${\bf F}_\ell^*$.
Then 
\[
\chi_1[G_\ell]\times \chi_2[G_\ell] \quad\text{equals}\quad 
C_2\times C_m \quad\text{or}\quad C_m\times C_2
\]
with $m>1$ dividing $\ell-1$. 
In the case $N_\ell=\ker\chi_1$ one has 
$\chi_1[N_\ell]\times \chi_2[N_\ell]\subset 1\times C_m$, and
in the case $N_\ell=\ker\chi_2$ has
$\chi_1[N_\ell]\times \chi_2[N_\ell]\subset C_m\times 1$.
We will count the number of pairwise non-conjugate pairs 
$(G_\ell, N_\ell)$ for each divisor $m>1$ of $\ell-1$.

If $m>1$ is odd, the inclusion \eqref{eqn:Goursatinclusion}
is an equality by the Goursat lemma, as $C_2$ and $C_m$ 
have no non-trivial common quotients.
In this case there are 3 pairwise non-conjugate pairs 
$(G_\ell, N_\ell)$. 
Indeed, in the `diagonal case' $G_\ell\iso \pi[G_\ell]$ when 
$T_\ell$ is trivial, the group $G_\ell$ stabilizes
two different lines on which it acts via $\chi_1$ and $\chi_2$.
In this case the two choices for $N_\ell=\ker\chi_i$ are
conjugate under a linear map that swaps these two lines, and
$(G_\ell, N_\ell)$ is unique up to conjugation in $\GL(V_\ell)$.
In the other case where $T_\ell\ne 1$ is non-trivial, 
$G_\ell$ stabilizes a unique line
corresponding to $\chi_1$, so the two choices for 
$N_\ell=\ker\chi_i$ give rise to non-conjugate pairs.

If $m>2$ is even, the groups $C_2\times C_m$ and $C_m\times C_2$ contain a 
unique index 2 subgroup that maps surjectively to each of the components. 
It corresponds to the Goursat isomorphism
$C_m/C_{m/2}\iso C_2$.
This gives 2 possibilities (of different order) for $\pi[G_\ell]$ 
inside each of $C_2\times C_m$ and $C_m\times C_2$.
Just as in the previous case of $m$ odd, these 4 possibilities
for $\pi[G_\ell]$ give rise to 4 `non-diagonal' possibilities for
$G_\ell$ with corresponding subgroups 
$N_\ell=\ker\chi_i\subset G_\ell$
that are pairwise non-conjugate because the 
non-trivial transvection subgroup $T_\ell\subset G_\ell$
creates an asymmetry between the characters $\chi_1$ and $\chi_2$ on $G_\ell$.
In the diagonal case, the swapping of the invariant lines 
corresponding to $\chi_1$ and $\chi_2$ reduces the number of
conjugacy classes in $\GL(V_\ell)$ to 2.

In the case $m=2$ excluded above, a small adaptation 
is necessary stemming from the fact that for $i(G)=2$,
the group $\pi[G_\ell]=\langle-1\rangle$ is so small that
in the non-diagonal case,
both choices $N_\ell=\ker\chi_i$ yield $N_\ell=1$, so they give rise to 
the same group up to conjugacy.
The table below summarizes the situation.
\begin{table}[h]
\centering
\renewcommand{\arraystretch}{1.5}   
\setlength{\tabcolsep}{10pt}       
\begin{tabular}{cccc|c}
\hline
$I(G)$ & $m$ & $T_\ell \neq 1$ & $T_\ell = 1$ & Total \\
\hline
$1$ & $m>1$ odd & $2$ & $1$ & $3$ \\
\hline
$1$ & \multirow{2}{*}{$m=2$} & $2$ & $1$ & $3$ \\
$2$ &                        & $1$ & $1$ & $2$ \\
\hline
$1$ & \multirow{2}{*}{$m>2$ even} & $2$ & $1$ & $3$ \\
$2$ &                             & $2$ & $1$ & $3$ \\
\hline
\end{tabular}
\end{table}

\noindent
To finish the proof, we add the values for the divisors 
$m>1$ of $\ell-1$ given in the table to find that up to conjugation,
the total number of special subgroups of level $2\ell$ equals
\smallskip\noindent
\[
\textstyle
3\sum_{m|\ell-1 \text{ odd}, m>1} 1 + 5 +
6\sum_{m|\ell-1 \text{ even}, m>2} 1 = 
\qquad\qquad\qquad\qquad
\]
\[
\textstyle
3\sum_{m|\ell -1 \text{ odd}} 1 - 4 +
6\sum_{m|\ell -1 \text{ even}} 1 
= 3\,\sigma_0(\ell-1)+3\,\sigma_0(\frac{\ell-1}{2}) - 4.
\]
\smallskip\noindent
This finishes the proof of Theorem \ref{thm:countingspecials}.\qed

\section{Elliptic curves over \texorpdfstring{$\mathbf Q$}{Q} giving rise to special groups}
\label{S:ellipticcurves}
\noindent
So far, our discussion has been purely group-theoretic, focusing 
on the classification and enumeration of special subgroups 
\[
    G_{2\ell} \subset \GL(V_2) \times \GL(V_\ell)  
\]
of level $2\ell$. As we indicated after Definition \ref{def:special_group},
the interest of these groups for us stems from the fact that 
the Galois group over a number field $K$ of the $2\ell$-division 
field $K_{2\ell}=K(E[2\ell](\overline K))$ of an elliptic curve $E/K$ 
is naturally a subgroup of $\GL(V_2) \times \GL(V_\ell)$: 
it acts linearly on the groups
$V_2$ and $V_\ell$ of 2- and $\ell$-torsion points of $E$.

Just as for the inverse Galois problem, there is the trivial
observation that every special subgroup of level $2\ell$ occurs
as a Galois group $\Gal(K_{2\ell}/K)$.
It suffices to take a number field $F$ and an elliptic curve 
$E/F$ for which the the $2\ell$-division field has the generic 
Galois group $\GL(V_2) \times \GL(V_\ell)$ over $F$, 
and then consider $E$
over the field $K$ that is left invariant by the 
subgroup that one wants to realize.

A more interesting question is: which special subgroups can be realized
for $E$ over $K=\Que$ as $\Gal(K_{2\ell}/K)$?
Their number is finite because of the following theorem.
\begin{theorem} \label{thm:special_elliptic_curves}
Let $\ell$ be an odd prime and $E/\Que$ an elliptic curve for which the Galois group over $\Que$ of the $2\ell$-division field is special. Then we have $\ell \in \{3,5,7\}$, and
$E$ is $\mathbf{Q}$-isogenous to the quadratic twist of an elliptic curve with a rational $\ell$-torsion point. 
\end{theorem}

\begin{proof}
Let $E/\Que$ be as in the Theorem.
By Theorem \ref{thm:2l_entanglement}, there is a $\Que$-rational 
subgroup $L\subset V_\ell=E[\ell](\overline\Que)$ of order 
$\ell$ such that the special Galois group
$G_{2\ell}$ of the $2\ell$-division field of $E$ over $\Que$
acts on $L$ by a character $\chi_1: G_{2\ell}\to {\bf F}_\ell^*$
and on $V_\ell/L$ by a character $\chi_2: G_{2\ell}\to {\bf F}_\ell^*$.
Moreover, one of these characters is quadratic.

If $\chi_1$ is quadratic, we can twist $E$ 
by $\chi_1$ to obtain a curve $E'/\Que$ 
with a $\Que$-rational point of order $\ell$. 
If not, the character $\chi_2$ is quadratic, and we can twist
$E$ by $\chi_2$ to obtain a curve $E'/\Que$ with a $\Que$-rational 
subgroup $L$ of order $\ell$ on which Galois acts via 
$\chi_1\chi_2$ and for which the $\Que$-isogenous curve
$E'/L$ has a $\Que$-rational point of order~$\ell$. 
As the only odd primes $\ell$ for which elliptic curves over $\Que$
can have points of order $\ell$ are 3, 5 and 7,
we have proved the Theorem.
\end{proof}

\noindent
The proof of Theorem \ref{thm:special_elliptic_curves} shows that
there is an easy way to construct elliptic curves over $\Que$ with 
special Galois group $G_{2\ell}$ for $\ell \in \{3,5,7\}$.
One simply performs the steps in the proof in reverse order,
starting from an elliptic curve $E_\ell/\mathbf{Q}$ with a rational 
$\ell$-torsion point $P$ and maximal 2-division field $K_2$. 
Such curves $E_\ell/\Que$ may be obtained by appropriately
specializing $t\in \mathbf{Q}$ in the Tate normal form \cite{Kubert}
\[
E_{\ell,t}:\quad y^2+axy+by=x^3+cx^2,
\]
where \(P=(0,0)\) has order \(\ell\), and
\[
(a,b,c)=
\begin{cases}
(1,t^{-1},0), &\text{for\ } \ell=3, \ t\neq 0,27; \\
(1-t,-t,-t), &\text{for\ }  \ell=5, \ t\neq 0; \\
(1+t-t^2,t^2-t^3,t^2-t^3), &\text{for\ }  \ell=7, \ t\neq 0,1.
\end{cases}
\]
The Galois group of the $\ell$-division field of the curve 
$E_\ell$ obtained in this way has a mod-$\ell$ representation 
of the form 
\[
\sigma \mapsto 
\begin{pmatrix}
1& *\\0&\omega_{\ell}(\sigma) 
\end{pmatrix}
\in\GL(V_\ell)
\]
in any basis for $V_\ell=E_\ell[\ell](\overline{\mathbf{Q}})$ for which 
the first vector spans the line through 0 and~$P$.
Here, $\omega_\ell: G_\Que \to \mathbf{F}_\ell^*$ denotes 
the mod-$\ell$ cyclotomic character; for elliptic curves 
over $\Que$, it is surjective and of order $\ell-1$.
The entry $*$ is non-canonical and depends on the choice of a second basis vector.

For the mod-2 representation of $E_\ell$, the character 
$\varepsilon: G_\mathbf{Q} 
    \twoheadrightarrow \{\pm 1\} \subset \mathbf{F}_\ell^*$ 
giving the sign of the permutation of the points of order 2
is different from $\omega_\ell$ for $\ell>3$, and also for
$\ell=3$ in cases where the 2-division field $K_2$ of $E_\ell$
does not contain $\mathbf{Q}(\zeta_3)$.
In this generic case, twisting $E_\ell$ by $\varepsilon$ yields 
an elliptic curve $E/\mathbf{Q}$ with the same mod-2-representation as $E_\ell$.
The Galois group $G_\ell$ of the mod-$\ell$ representation
\[
\sigma \mapsto 
\begin{pmatrix}
\varepsilon(\sigma) & *\\0&\varepsilon(\sigma) \omega_{\ell}(\sigma) 
\end{pmatrix}
\in\GL(V_\ell)
\]
of $E$ is a Borel subgroup of $\GL(V_\ell)$ with fundamental character
$\chi_1=\varepsilon$.
In this case, the Galois group $G_{2\ell}$ of the $2\ell$-division field of $E$ over $\Que$
is the special subgroup of level  $2\ell$ described by \eqref{eqn:def-of-G}, with $i=1$.

\begin{example}
The curve $E_3:y^2=x^3+x^2+2x+1$ with LMFDB label \href{https://www.lmfdb.org/EllipticCurve/Q/92/b/2}{92.b2} 
has discriminant $\Delta=-2^4 \cdot 23$ and a rational 3-torsion point $P=(0,1)$.
Its 2-division field $K_2$ has degree 6 over $\mathbf{Q}$
and quadratic subfield $\Que(\sqrt{\Delta})=\Que(\sqrt{-23})$. 

The Galois group $\Gal(\mathbf{Q}(E_3[3](\overline{\mathbf{Q}}))/\mathbf{Q})$ is contained in 
the Borel subgroup relative to the line $L={\bf F}_3\cdot P$, with fundamental characters
$\chi_1=1$ and $\chi_2=\omega_3$. 

As the Dirichlet character $\alpha$ corresponding to $\Que(\sqrt{-23})$ is different from $\omega_3$,
the curve $E=E_3^\alpha$ obtained by twisting $E_3$ with $\alpha$ 
has 6-division field $K_6$ with 
\begin{equation}\label{eqn:G6}
G_6= \Gal(K_6/\mathbf{Q})\ \mathop{\subset}\limits_2 \ G_2\times G_3=\GL(V_2) \times B
\end{equation}
a special Galois group of level 6.
Here $B=G_3=\Gal(K_3/\Que)\subset \GL(V_3)$ is the Borel subgroup with 
$\chi_1=\alpha$ and $\chi_2=\alpha \omega_3$.
The entanglement in \eqref{eqn:def-of-G} is given by $\chi_1=\alpha$ corresponding to
$\Que(\sqrt{-23})=K_2\cap K_3$.

The curve $E$ has twisted Weierstrass equation 
$-23y^2= x^3+x^2+2x+1$,

and the following table gives the orders of the groups $E(\mathbf{F}_p)$ for primes $p\leq 97$.
Note that as expected, these orders can be coprime to 2 and coprime to 3, but not to 6.

\begin{table}[ht]
\centering
\renewcommand{\arraystretch}{1.3}
\setlength{\tabcolsep}{5pt}
\begin{tabular}{c|*{8}{c}}
$p$ & $3$ & $5$ & $7$ & $11$ & $13$ & $17$ & $19$ & $29$ \\
\hline
$\#E(\mathbf{F}_p)$ & 3 & $2\cdot 3$ & $2\cdot 5$ & $2^2\cdot 3$ & $3\cdot 5$ & $2^2\cdot 3$ & $2\cdot 11$ & $3\cdot 11$ \\
\noalign{\vskip 8pt}
$p$ & $31$ & $37$ & $41$ & $43$ & $47$ & $53$ & $59$ & $61$ \\
\hline
$\#E(\mathbf{F}_p)$  & $3^3$ & $2\cdot 23$ & $3\cdot 13$ & $2^2\cdot 13$ & $3\cdot 13$ & $2^2\cdot 3\cdot 5$ & $2^3\cdot 3^2$ & $2^2\cdot 19$ \\
\noalign{\vskip 8pt}
$p$ & $67$ & $71$ & $73$ & $79$ & $83$ & $89$ & $97$ &  \\
\hline
$\#E(\mathbf{F}_p)$ & $2^2\cdot 19$ & $3\cdot 29$ & $3^4$ & $2\cdot 5\cdot 7$ & $2\cdot 3^2\cdot 5$ & $2\cdot 3^2\cdot 5$ & $2^3\cdot 11$ &  \\
\end{tabular}

\end{table}
\end{example}
\noindent
For the curve $E_3$ and the line $L={\bf F}_3\cdot P$ generated by the rational 3-torsion point~$P$ we can compute
the 3-isogenous curve $\widetilde E=E_3/L$ as in \cite{Velu}.
It has Weierstrass equation $y^2=x^3+x^2-18x-43$ and LMFDB label
\href{https://www.lmfdb.org/EllipticCurve/Q/92/b/1}{92.b1}.
The Galois group over $\Que$ of its 3-division field $\widetilde K_3$ is 
a Borel subgroup with fundamental characters
$\chi_1=\omega_3$ and $\chi_2=1$, and its 2-division field is the same as for $E_3$.
The twist $\widetilde E^\alpha$ with twisted Weierstrass equation $-23y^2=x^3+x^2-18x-43$
has a special Galois group of level $6$ that is described as in \eqref{eqn:G6}.
This time the Borel subgroup $B$ has $\chi_1=\alpha\omega_3$ and $\chi_2=\alpha$, and the entanglement 
in \eqref{eqn:def-of-G} is via $\chi_2$.
The curve $\widetilde E^\alpha$ is 3-isogenous to the curve $E=E_3^\alpha$ above, so the orders of their point groups
over ${\bf F}_p$ for primes of good reduction $p\nmid 2\cdot 23$ coincide.

\end{document}